# The Strong Macdonald Conjecture

S. Fishel, I. Grojnowski and C. Teleman

**Introduction**

This note presents our proof of the *strong Macdonald constant term conjecture* of Hanlon [H] and Feigin [F1] for reductive Lie algebras. It is part of a larger paper, on the Dolbeault cohomology of flag varieties of loop groups and combinatorial applications to certain hypergeometric summation formulas; results for $\mathfrak{sl}_2$ are described in [T2]. (Some of our identities were independently discovered, in more general form, by Macdonald [M]). However, as the proof of the conjecture is of interest independent of the rest of the paper, and given the delays in completing other parts of our project, we have decided to make this self-contained section available for distribution. Understandably, we refrain from extensive background discussions, which will find their place in the complete paper.

**2. The Strong Macdonald Conjecture**

*2.1. Notations*
Let $\mathfrak{g}$ be the Lie algebra of a complex reductive Lie group $G$ of rank $\ell$, with *exponents* $m_1,\ldots,m_\ell$. Recall, when $\mathfrak{g} = \mathfrak{gl}_n$, that $\ell = n$ and $(m_1,\ldots,m_n) = (0,\ldots,n-1)$. In general, two related definitions for the exponents come from the following statements:
  • the algebra $(S\mathfrak{g})^\mathfrak{g}$ of invariant polynomials on $\mathfrak{g}^*$, under the co-adjoint action, is a free symmetric algebra generated in degrees $m_1 + 1,\ldots,m_\ell + 1$;
  • the subalgebra $(\Lambda\mathfrak{g})^\mathfrak{g}$ of ad-invariants in the full exterior algebra of $\mathfrak{g}$ is a free exterior algebra generated in degrees $2m_1 + 1,\ldots,2m_\ell + 1$.

Call $\mathfrak{g}[z]$ the Lie algebra of $\mathfrak{g}$-valued polynomials in $z$, $\mathfrak{g}[z]/z^n$ the truncated polynomial Lie algebra, and $\mathfrak{g}[z,s] := \mathfrak{g} \otimes \mathbb{C}[z,s]$ the *graded* Lie algebra based on the free skew–commutative algebra $\mathbb{C}[z,s]$, with the odd variable $s$ of degree 1. The third one is a semi-direct product $\mathfrak{g}[z] \ltimes s\mathfrak{g}[z]$, with the adjoint action of the first factor on the second. We view them all as Lie algebras over $\mathbb{C}$.

Choose, for later use, a compact form of $G$, and, when $\mathfrak{g}$ is semi-simple, a basis of self-adjoint elements $\xi_a$ in $\mathfrak{g}$, orthonormal with respect to the Killing form. Call, for $m \geq 0$, $\psi^a(-m)$ and $\sigma^a(-m)$ the elements of $\Lambda^1\mathfrak{g}[z]^*$ and of $S^1\mathfrak{g}[z]^*$ dual to the basis $z^m \cdot \xi_a$ of $\mathfrak{g}[z]$.

For any complex Lie algebra $\mathfrak{L}$, $H^*(\mathfrak{L})$ will denote its *cohomology with complex coefficients*. Conceptually, it is $\mathrm{Ext}^*_\mathfrak{L}(\mathbb{C};\mathbb{C})$ in the category of complex $\mathfrak{L}$-representations ($\mathfrak{L}$-modules), but we shall define it via its *Koszul complex* (2.6). $H^*(\mathfrak{L};V)$ shall stand for *cohomology with coefficients* in an $\mathfrak{L}$-*module V*; it is $\mathrm{Ext}^*_\mathfrak{L}(\mathbb{C};V)$. If $\mathfrak{g} \subseteq \mathfrak{L}$ is a subalgebra, we call $H^*(\mathfrak{L},\mathfrak{g};V)$ the *relative* Lie algebra cohomology. If $\mathfrak{g}$ is reductive and $\mathfrak{L}$ (via ad) and $V$ (by restriction) are semi-simple $\mathfrak{g}$-modules, this is again $\mathrm{Ext}^*_\mathfrak{L}(\mathbb{C};V)$, but in the category of $\mathfrak{g}$-semi-simple $\mathfrak{L}$-modules.

*2.2. Background*
The "strong Macdonald" conjectures describe the cohomologies of the Lie algebra $\mathfrak{g}[z]/z^n$ and of the graded Lie algebra $\mathfrak{g}[z,s]$. We call the two versions of the conjecture the *truncated algebra* form (2.3.1) and the *super-algebra* form (2.3.3). The truncated form is due to Hanlon [H], who proved it for $\mathfrak{gl}_n$. It may have been independently known to Feigin [F1], who is responsible for the super-algebra version, and for relating the two [F2].

Feigin also outlined an argument for the super-algebra version; but we are unsure whether it can be carried out as indicated[1]. While we cannot fill in the gap, we do confirm the conjectures by a different

---
[1] One particular step, the lemma on p. 93 of [F2], is incorrect; it fails for absolute cohomology when $Q = \partial/\partial\xi$, and nothing in the suggested argument seems to account for that

route: we compute the cohomology of $\mathfrak{g}[z,s]$ directly, by finding the *harmonic cocycles* in the Koszul complex, in a suitable metric. That being accomplished, Feigin's argument easily recovers the cohomology of the truncated Lie algebra.

The success of our Laplacian approach comes from the specific metric used on the Koszul complex, which originates in the Kähler geometry of the *flag variety* of the loop group of *G*. The latter is responsible for an identity between two different Laplacians, far from obvious in explicit form, which implies that the harmonic cocycles form a *subalgebra*, and allows their computation. (This is closely related to results of [FGZ]). We do not know of a computation using the more obvious (Killing-like) metric; its harmonic cocycles are not closed under multiplication.

*2.3. Statements*

The following affirms Hanlon's conjecture for reductive $\mathfrak{g}$.

**(2.3.1) Theorem.** $H^*\left(\mathfrak{g}[z]/z^n\right)$ *is a free exterior algebra on* $n \cdot \ell$ *generators, with n generators in cohomology degree* $2m+1$ *for each exponent* $m = m_1,\ldots,m_\ell$. *The z-weights of these generators are the negatives of:* $0, mn+1, mn+2,\ldots, mn+n-1$. ✠

*(2.3.2) Remarks.* (i) Reductivity of $\mathfrak{g}$ and the semi-direct decomposition $\mathfrak{g}[z]/z^n = \mathfrak{g} \ltimes z\mathfrak{g}[z]/z^n$ lead to a *natural* isomorphism $H^*\left(\mathfrak{g}[z]/z^n\right) = H^*(\mathfrak{g}) \otimes H^*\left(\mathfrak{g}[z]/z^n, \mathfrak{g}\right)$. The first factor accounts for the generators of *z*-weight 0.
(ii) Ignoring *z*-weights leads to an *abstract* ring isomorphism $H^*\left(\mathfrak{g}[z]/z^n\right) \cong H^*(\mathfrak{g})^{\otimes n}$.
(iii) The dimensional lower bound $H^*\left(\mathfrak{g}[z]/z^n\right) \geq H^*(\mathfrak{g})^{\otimes n}$ holds for *any* Lie algebra $\mathfrak{g}$. Namely, $\mathfrak{g}[z]/z^n$ is a degeneration of $\mathfrak{g}[z]/(z^n - \varepsilon)$ as $\varepsilon = 0$. When $\varepsilon \neq 0$, the latter is isomorphic to $\mathfrak{g}^{\oplus n}$, whose cohomology is $H^*(\mathfrak{g})^{\otimes n}$; and the cohomology can only jump up under degeneration. However, this argument says little about the *ring* structure.

A conceptual formulation of theorem (2.3.1) was suggested independently by Feigin and Loday. Given a skew-commutative algebra *A* and any Lie algebra $\mathfrak{g}$, an invariant polynomial $\Phi$ of degree $(m+1)$ on $\mathfrak{g}$ determines a linear map from the dual of $HC_n^{(m)}(A)$, the *m*th *Adams component* of the *n*th *cyclic homology* of *A*, to $H^{n+1}(\mathfrak{g} \otimes A)$. (See (2.4.3) below for the case of interest here, or [T2], (2.2) in general). When $\mathfrak{g}$ is reductive, Loday suggested that these maps might be injective, and that $H^*(\mathfrak{g} \otimes A)$ might be freely generated by their images, as $\Phi$ ranges over a set of generators of the ring of invariant polynomials. (The Adams degree *m* will then range over the exponents $m_1,\ldots,m_\ell$). Thus, for $A = \mathbb{C}$, $HC_n^{(m)} = 0$ for $n \neq 2m$, while $HC_{2m}^{(m)} = \mathbb{C}$; and we recover the well-known description of $H^*(\mathfrak{g})$. For $\mathfrak{g} = \mathfrak{gl}_\infty$ and any associative, unital, graded *A*, this is the theorem of Loday-Quillen [LQ] and Tsygan [Ts]. It emerges from its proof that (2.3.1) affirms Loday's conjecture for $\mathbb{C}[z]/z^n$, while (2.3.3) below does the same for the graded algebra $\mathbb{C}[z,s]$. (The conjecture fails in general, and no characterization is known of the algebras *A* for which it holds; see [T2] for more details).

Describing the cohomology of the infinite-dimensional Lie algebra $\mathfrak{g}[z,s]$ requires some preparation. We start with its *homology*, which has a natural *co-algebra* structure. (Our general discomfort with co-algebras will soon make us switch back to cohomology). We can factor $H_*\left(\mathfrak{g}[z,s]\right)$ naturally as $H_*(\mathfrak{g}) \otimes H_*(\mathfrak{g}[z,s], \mathfrak{g})$, for the reason described in (2.3.2.i); the first factor behaves rather differently from the rest, and is best set aside.

**(2.3.3) Theorem.** $H_*\left(\mathfrak{g}[z,s], \mathfrak{g}\right)$ *is isomorphic to the free graded co-commutative co-algebra whose space of primitives is the direct sum of* $\mathbb{C}[z] \cdot s^{\otimes(m+1)}$, *in total degree* $2m+2$, *and* $\mathbb{C}[z]dz \cdot s^{\otimes m}$, *in total degree* $2m+1$, *as m ranges over the exponents* $m_1,\ldots,m_\ell$. *The isomorphism respects* $(z,s)$-*weights*. ✠

*(2.3.4) Remarks.* (i) As a vector space, a free co-commutative co-algebra is isomorphic to the graded symmetric algebra on its primitives; but there is no *a priori* algebra structure on homology.
(ii) The *total degree* includes the degree of $s$. The degrees of the cycles as multi-linear forms are $m+1$ in both cases: see their description (2.4.3) below.

The description in (2.3.3) is not quite canonical, although it does correctly indicate the infinitesimal automorphisms induced by polynomial vector fields in $z$. If $P^{(m)}$ is the $m$th degree part of the space of primitives in the quotient co-algebra $S\mathfrak{g}/[\mathfrak{g}, S\mathfrak{g}]$, canonical descriptions of our primitives are

(2.3.5') $\qquad \bigoplus_m P^{(m+1)} \otimes \mathbb{C}[z] \cdot s(ds)^m,$

(2.3.5'') $\qquad \bigoplus_m P^{(m+1)} \otimes \left( \mathbb{C}[z] \cdot (ds)^m + \mathbb{C}[z]dz \cdot s(ds)^{m-1} \right) \big/ d\left( \mathbb{C}[z] \cdot s(ds)^{m-1} \right).$

The right factors are the cyclic homology components $HC_{2m+1}^{(m)}$ and $HC_{2m}^{(m)}$ of the *non-unital algebra* $\mathbb{C}[z,s] \ominus \mathbb{C}$. (The last factor, $HC_{2m}^{(m)}$, is identifiable with $\mathbb{C}[z]dz \cdot s(ds)^{m-1}$, for $m \neq 0$, and with $\mathbb{C}[z]/\mathbb{C}$ if $m = 0$.) This description is compatible with the action of *super-vector fields* in $z$ and $s$; see (2.5.5) below.

*2.4. Restatement without super-algebras*

There is a natural isomorphism [K] between $H_*(\mathfrak{L}; \Lambda V)$ and the homology of the semi-direct product Lie algebra $\mathfrak{L} \ltimes V$, with trivial Lie bracket on $V$. Its graded version, applied to $\mathfrak{L} = \mathfrak{g}[z]$ and the (odd) vector space $V = s\mathfrak{g}[z]$, is the following equality:

(2.4.1) $\qquad H_n(\mathfrak{g}[z,s], \mathfrak{g}) = \bigoplus_{p+q=n} H_{q-p}(\mathfrak{g}[z], \mathfrak{g}; S^p(s\mathfrak{g}[z]));$

the generators from $s\mathfrak{g}[z]$ carry homology degree 2, combining their degree as a form with that of $s$.

**(2.4.2) Theorem.** $H_*(\mathfrak{g}[z], \mathfrak{g}; S(s\mathfrak{g}[z]))$ *is the free graded co-commutative co-algebra with primitive elements* $\mathbb{C}[z] \cdot s^{\otimes(m+1)}$, *in degree 0, and primitive elements* $\mathbb{C}[z]dz \cdot s^{\otimes m}$ *in degree 1, as $m$ ranges over the exponents $m_1, \ldots, m_\ell$. This preserves the $z$- and $s$-weights.* ✠

While $H^*(\mathfrak{g}[z,s], \mathfrak{g})$ is obtained by straight duality, infinite-dimensionality makes it a bit awkward, and we opt for restricted duality, defined from the direct sum of the $(s,z)$-weight spaces in the dual of the Koszul complex of (2.6) below. These weight spaces are finite-dimensional and are preserved by the Koszul differential. (A single grading with positive weights on $s$ and $z$ would have the same effect). The cohomology of this restricted complex, called *restricted Lie algebra cohomology* $H^*_{res}(\mathfrak{g}[z,s], \mathfrak{g})$, is the sum of weight spaces in the full dual $H^*$ of (2.4.1).

**(2.4.3) Theorem.** $H^*_{res}(\mathfrak{g}[z], \mathfrak{g}; S\mathfrak{g}[z]^*_{res})$ *is isomorphic to the free graded commutative algebra generated by the restricted duals of $\bigoplus_m P^{(m+1)} \otimes \mathbb{C}[z]$ and $\bigoplus_m P^{(m+1)} \otimes \mathbb{C}[z]dz$, in degrees 0 and 1, respectively. Specifically, an invariant linear map $\Phi: S^{m+1}\mathfrak{g} \to \mathbb{C}$ determines linear maps*

$$S_\Phi: S^{m+1}\mathfrak{g}[z] \to \mathbb{C}[z], \qquad \sigma_0 \cdot \sigma_1 \cdot \ldots \cdot \sigma_m \mapsto \Phi(\sigma_0(z), \sigma_1(z), \ldots, \sigma_m(z)),$$

$$E_\Phi: \Lambda^1(\mathfrak{g}[z]/\mathfrak{g}) \otimes S^m\mathfrak{g}[z] \to \mathbb{C}[z]dz, \quad \psi \otimes \sigma_1 \cdot \ldots \cdot \sigma_m \mapsto \Phi(d\psi(z), \sigma_1(z), \ldots, \sigma_m(z)).$$

*Therein, the coefficients $S_\Phi(-n)$, $E_\Phi(-n)$ of $z^n$, resp. $z^{n-1}dz$ give restricted 0- and 1-cocycles, and $H^*_{res}$ is freely generated by these, if $\Phi$ ranges over a generating set of invariant polynomials on $\mathfrak{g}$.* ✠

To illustrate, here are the cocycles associated to the Killing form on $\mathfrak{g}$; see (2.1) for notations.

$$S(-n) = \sum_{\substack{1 \leq a \leq \dim \mathfrak{g} \\ 0 \leq p \leq n}} \sigma^a(-p)\sigma^a(p-n), \qquad E(-n) = \sum_{\substack{1 \leq a \leq \dim \mathfrak{g} \\ 0 < p \leq n}} p\psi^a(-p)\sigma^a(p-n).$$

Let now $\mathfrak{B} \subset \mathfrak{g}[z]$ be an *Iwahori subalgebra*, inverse image of a Borel subalgebra $\mathfrak{b} \subset \mathfrak{g}$ under the evaluation $\mathfrak{g}[z] \to \mathfrak{g}$ at $z = 0$. Let also $\mathfrak{h} \subset \mathfrak{b}$ be a Cartan subalgebra. Closely related to (2.4.2) is a description of the cohomology of the super-algebra pair $(\mathfrak{B}[s], \mathfrak{h})$.

Note that he cocycles $S_\Phi(0)$ generate a copy of $(S(s\mathfrak{g})^*)^\mathfrak{g}$ within $H^*_{res}(\mathfrak{g}[z,s], \mathfrak{g})$. A similar ring monomorphism $S(s\mathfrak{h})^* \to H^*_{res}(\mathfrak{B}[s], \mathfrak{h})$ arises by identifying $\mathfrak{h}^*$ with the $\mathfrak{B}$-invariants in $\mathfrak{B}^*$. Recall, finally, that $(S\mathfrak{g}^*)^\mathfrak{g}$ sits within $S\mathfrak{h}^*$ as the subring of Weyl invariants.

**(2.4.4) Theorem.** $H^*_{res}(\mathfrak{B}[s], \mathfrak{h}) = H^*_{res}(\mathfrak{g}[z,s], \mathfrak{g}) \otimes_{(S(s\mathfrak{g})^*)^\mathfrak{g}} S(s\mathfrak{h})^*$, *naturally.* ✠

In other words, when going from $\mathfrak{g}[z]$ to $\mathfrak{B}$, the factor $(S\mathfrak{g}^*)^\mathfrak{g}$ of zero-modes is replaced with $S\mathfrak{h}^*$.

*2.5. Proof of the truncated version*

Assuming the super-algebra result (2.4.3), Feigin's constructions in [F2] provide a proof of the truncated version (2.3.1) of the conjecture; their shadow is the specialization $t = q^n$ in the combinatorial literature (in our case, this is $s = z^n$).

The resolution of $\mathfrak{g}[z]/z^n$ by the *differential graded Lie algebra* $\mathfrak{L}_\bullet$, with differential $\partial s = z^n$

(2.5.1) $\qquad (\mathfrak{L}_1 := s\mathfrak{g}[z]) \xrightarrow{\partial : s \mapsto z^n} (\mathfrak{L}_0 := \mathfrak{g}[z]) \to \mathfrak{g}[z]/z^n$,

identifies $H^*(\mathfrak{g}[z]/z^n)$ with the *hypercohomology* of $\mathfrak{L}_\bullet$, and $H^*(\mathfrak{g}[z]/z^n, \mathfrak{g})$ with that of the pair $(\mathfrak{L}_\bullet, \mathfrak{g})$. Recall that hypercohomology is computed by a *double complex*, where the Koszul differential is supplemented by a differential induced by $\partial$. This leads to a convergent spectral sequence, with

(2.5.2) $\qquad E_1^{p,q} = H^{q-p}_{res}\left(\mathfrak{g}[z], \mathfrak{g}; S^p(s\mathfrak{g}[z])^*_{res}\right) \Rightarrow H^{p+q}\left(\mathfrak{g}[z]/z^n, \mathfrak{g}\right)$.

The $E_1^{p,q}$ term arises by ignoring $\partial$, and is the portion of $H^{p+q}_{res}(\mathfrak{g}[z,s], \mathfrak{g})$ with $s$-weight $(-p)$, cf. (2.4.1). If we assign weight 1 to $z$ and weight $n$ to $s$, then $\mathfrak{L}_\bullet$ carries this additional $z$-grading, preserved by $\partial$ and hence by the spectral sequence.

**(2.5.3) Lemma.** *Let $n > 0$. $E_2^{p,q}$ is the free skew-commutative algebra generated by the dual of the sum of vector spaces $s^{\otimes m} \mathbb{C}[z]dz/d(z^n \cdot \mathbb{C}[z])$, placed in bi-degrees $(p,q) = (m, m+1)$, as $m$ ranges over $m_1, \ldots, m_\ell$. The $z$-weight of $s$ is $n$.*

*Proof of Theorem (2.3.1).* All differentials, as derivations, are determined by their action on generators; the placement of the generators of $E_2$ in (2.5.3) disallows higher differentials, so the spectral sequence collapses here[2]. Therefore, $E_2 = E_\infty$ is the associated graded ring for a filtration on $H^*(\mathfrak{g}[z]/z^n, \mathfrak{g})$, which is compatible with the $z$-grading. Freedom of the associated graded ring forces $H^*(\mathfrak{g}[z]/z^n, \mathfrak{g})$ to be isomorphic to the same, and we get the desired description, in view of (2.3.2.i). ✠

*Proof of Lemma (2.5.3).* The description of the generating cocycles $E_\Phi$ and $S_\Phi$ of $E_1$ in (2.4.3) allow us to compute $\delta_1$. The $S_\Phi$ have nowhere to go, but for $E_\Phi : \Lambda^1 \otimes S^p \to \mathbb{C}[z]dz$, we get

(2.5.4)
$$(\delta_1 E_\Phi)(\sigma_0 \cdot \ldots \cdot \sigma_m) = E_\Phi(\partial(\sigma_0 \cdot \ldots \cdot \sigma_m)) = \sum_k E_\Phi(z^n \sigma_k \otimes \sigma_0 \cdot \ldots \cdot \hat{\sigma}_k \cdot \ldots \cdot \sigma_m) =$$
$$= \sum_k \Phi(\sigma_0 \cdot \ldots \cdot d(z^n \sigma_k) \cdot \ldots \cdot \sigma_m) =$$
$$= (m+1)n \cdot z^{n-1} dz \cdot \Phi(\sigma_0 \cdot \ldots \cdot \sigma_m) + z^n \cdot d\Phi(\sigma_0 \cdot \ldots \cdot \sigma_m) =$$
$$= \left((m+1)n \cdot z^{n-1} dz \cdot S_\Phi + z^n \cdot dS_\Phi\right)(\sigma_0 \cdot \ldots \cdot \sigma_m)$$

---

[2] For an alternative argument, cf. (2.3.2.iii): the second $E$ term already meets the dimensional lower bound.

and so $\delta_1$ is the transpose of the linear operator $(m+1)n \cdot z^{n-1}dz \wedge + z^n \cdot d$, from $\mathbb{C}[z]$ to $\mathbb{C}[z]dz$. This has no kernel for $n > 0$, and its cokernel is $\mathbb{C}[z]dz/d(z^n \cdot \mathbb{C}[z])$. ✠

*(2.5.5) Remarks.* (i) $\delta_1$ is induced by $\partial$, which on $\mathfrak{g}[z,s]$ is given by the super-vector field $z^n \partial/\partial s$. The latter acts naturally on the presentation (2.3.5'), (2.3.5") of the cycles:

(2.5.6) $$z^n \partial/\partial s : \mathbb{C}[z] \cdot s(ds)^m \to \frac{\mathbb{C}[z] \cdot (ds)^m + \mathbb{C}[z]dz \cdot s(ds)^{m-1}}{d(\mathbb{C}[z] \cdot s(ds)^{m-1})}.$$

Identifying the target space with $\mathbb{C}[z]dz \cdot s(ds)^{m-1}$ by projection, one checks that $z^k \cdot s(ds)^m$ maps to $(mn+n+k) \cdot z^{n+k-1}dz \cdot s(ds)^{m-1}$, which agrees with the transpose of $\delta_1$, as was claimed after (2.3.5).
(ii) If $n = 0$, the map in (2.5.6) is surjective, with 1-dimensional kernel; so $E_\infty^{p,q}$ now lives on the diagonal and equals $(S^p\mathfrak{g}^*)^\mathfrak{g}$. This is, in fact, a correct interpretation of $H^*(0,\mathfrak{g};\mathbb{C})$!

## 2.6. The Koszul complex

The *Lie algebra homology Koszul complex* [K] for a Lie algebra $\mathfrak{L}$ with coefficients in a module $V$ is $\Lambda\mathfrak{L} \otimes V$, homologically graded by the exterior degree, with differential

$$\delta(\lambda_1 \wedge \ldots \wedge \lambda_n \otimes v) = \sum_p (-1)^p \lambda_1 \wedge \ldots \wedge \hat{\lambda}_p \wedge \ldots \wedge \lambda_n \otimes \lambda_p(v)$$
$$+ \sum_{p<q} (-1)^{p+q} [\lambda_p, \lambda_q] \wedge \lambda_1 \wedge \ldots \wedge \hat{\lambda}_p \wedge \ldots \wedge \hat{\lambda}_q \wedge \ldots \wedge \lambda_n \otimes v$$

where the hats indicate missing factors. Its homology $H_*(\mathfrak{L};V)$ is the *Lie algebra homology of $\mathfrak{L}$ with coefficients in V*. If $\mathfrak{g} \subseteq \mathfrak{L}$ is a subalgebra, $\delta$ descends to the quotient $\Lambda(\mathfrak{L}/\mathfrak{g}) \otimes V/\mathfrak{g}.(\Lambda(\mathfrak{L}/\mathfrak{g}) \otimes V)$ of co-invariants under $\mathfrak{g}$. The full algebraic dual of its homology is the *relative cohomology* $H^*(\mathfrak{L},\mathfrak{g};V^*)$ of (2.1). It is computed by the dual complex, the *Koszul complex for Lie algebra cohomology of the pair $(\mathfrak{L},\mathfrak{g})$ with coefficients in $V^*$*. When $V^*$ has an algebra structure and $\mathfrak{L}$ acts by derivations, the Koszul complex is a differential graded algebra.

Let now $\mathfrak{g}$ be semi-simple, $\mathfrak{L} = \mathfrak{g}[z]$, $V = S\mathfrak{g}[z]$. For explicit work, we introduce the following operators, defined for any $m \in \mathbb{Z}$ by brutal truncation of the adjoint action of $\mathfrak{g}[z,z^{-1}]$, and extended to derivations of $\Lambda \otimes S := \Lambda(\mathfrak{g}[z]/\mathfrak{g})^*_{res} \otimes S\mathfrak{g}[z]^*_{res}$; the notations are those of (2.1).

(2.6.1) $$\mathrm{ad}_a(m) : \psi^b(n) \mapsto \begin{cases} \psi^{[a,b]}(m+n), & \text{if } m+n < 0 \\ 0, & \text{if } m+n \geq 0 \end{cases};$$

(2.6.2) $$R_a(m) : \sigma^b(n) \mapsto \begin{cases} \sigma^{[a,b]}(m+n), & \text{if } m+n \leq 0 \\ 0, & \text{if } m+n > 0 \end{cases};$$

$a,b$ range over $A := \{1, \ldots, \dim\mathfrak{g}\}$. The abusive notation $\psi^{[a,b]}(n)$ is self-explanatory. Finally, let

(2.6.3) $$\bar\partial = \sum_{a \in A; m > 0} \left(\psi^a(-m) \otimes R_a(m) + \psi^a(-m) \cdot \mathrm{ad}_a(m) \otimes 1/2\right),$$

where $\psi^a(-m)$ doubles notationally for the appropriate multiplication operator. (The notation for $\bar\partial$ originates in its geometric interpretation as a Dolbeault operator.)

**(2.6.4) Definition.** *The restricted Koszul complex $(C^\bullet, \bar\partial)$ for the pair $(\mathfrak{g}[z],\mathfrak{g})$ with coefficients in $S\mathfrak{g}[z]^*_{res}$ consists of the $\mathfrak{g}$-invariant skew polynomials in the $\psi^a(m)$ ($m < 0$) and the $\sigma^a(m)$ ($m \leq 0$), co-homologically graded by the $\psi$'s, with differential (2.6.3).*

*2.7. The metric and the Laplacian*

Define a hermitian metric on the Koszul complex by setting, if $m = n$ and $a = b$, $\langle \sigma^a(m) | \sigma^b(n) \rangle = 1$, $\langle \psi^a(m) | \psi^b(n) \rangle = -1/n$, and both products to zero otherwise, and then taking the multi-linear extension to $\Lambda \otimes S$. (Thus, $\sigma^a(m)^n / \sqrt{n!}$ has norm 1). The hermitian adjoints to (2.6.1) are the derivations defined by

(2.6.1*) $$\mathrm{ad}_a(m)^* \left( \psi^b(n) \right) = \frac{n-m}{n} \psi^{[a,b]}(n-m), \text{ or zero, if } n \geq m.$$

The $R$'s of (2.6.2) satisfy the simpler relation $R_a(m)^* = R_a(-m)$. The adjoint of (2.6.3) is

(2.6.3*) $$\bar{\partial}^* = \sum_{a \in A; m > 0} \left( \psi^a(-m)^* \otimes R_a(-m) + \mathrm{ad}_a(m)^* \cdot \psi^a(-m)^* \otimes 1/2 \right).$$

Finally, introduce the Laplacian $\bar{\Box} = (\bar{\partial} + \bar{\partial}^*)^2 = \bar{\partial}\bar{\partial}^* + \bar{\partial}^*\bar{\partial}$. It preserves the cohomology grading.

A (restricted) Koszul cocycle in the kernel of $\bar{\Box}$ is called *harmonic*. Noting that $\bar{\partial}$, $\bar{\partial}^*$ and $\bar{\Box}$ preserve the orthogonal decomposition into the finite-dimensional $(z,s)$-weight spaces, elementary linear algebra gives the following "Hodge decomposition":

**(2.7.1) Proposition.** *The map from harmonic cocycles $\mathcal{H}^k \subset C^k$ to their cohomology classes, via the decompositions $\ker \bar{\partial} = \mathrm{Im}\, \bar{\partial} \oplus \mathcal{H}^k$, $C^k = \mathrm{Im}\, \bar{\partial} \oplus \mathcal{H}^k \oplus \mathrm{Im}\, \bar{\partial}^*$, is a linear isomorphism.* ✠

We now set out to find $\mathcal{H}^k$, to which end we introduce the following adjoint pairs of operators:

(2.7.2)   $d_a(m): \sigma^b(n) \mapsto \psi^{[a,b]}(m+n)$, or zero, if $m+n \geq 0$,   $d_a(m)\psi^b(n) = 0$

(2.7.2*)   $d_a(m)^*: \psi^b(n) \mapsto -\sigma^{[a,b]}(n-m)/n$, or zero, if $n > m$,   $d_a(m)^* \sigma^b(n) = 0$

extended to (odd degree) derivations of $\Lambda \otimes S$. Finally, the following operators evidently preserve $C^\bullet$:

(2.7.3) $$D := \sum_{m > 0; a \in A} d_a(-m) d_a(-m)^*$$

(2.7.4) $$\Box := \sum_{a \in A; m > 0} \frac{1}{m} \left( R_a(-m) + \mathrm{ad}_a(-m) \right) \left( R_a(m) + \mathrm{ad}_a(-m)^* \right).$$

**(2.7.5) Theorem.** *On $C^\bullet$, $\bar{\Box} = \Box + D$. In particular, the harmonic forms are the joint kernel in $\Lambda \otimes S$ of the derivations $d_a(-m)^*$, as $a \in A$, $m > 0$, and $R_a(m) + \mathrm{ad}_a(-m)^*$, as $a \in A$, $m \geq 0$.*

We shall exploit this in the next subsection, where we identify the harmonic forms. Note, for now, that they must form a subalgebra, since they are cut out by derivations. The rest of this subsection is devoted to proving Theorem (2.7.5). We give two proofs.

*First proof.* Introduce yet another operator

(2.7.6) $$K := \sum_{a, b \in A; m > 0} \left( R_{[a,b]}(0) + \mathrm{ad}_{[a,b]}(0) \right) \cdot \psi^a(-m) \wedge \psi^b(-m)^*.$$

Note that the $\psi \wedge \psi^*$ term could equally well be written in first position, because

$$\sum_{a,b} \left[ \mathrm{ad}_{[a,b]}(0), \psi^a(-m) \cdot \psi^b(-m)^* \right] = \sum_{a,b} \left( \psi^{[[a,b],a]}(-m) \wedge \psi^b(-m)^* + \psi^a(-m) \wedge \psi^{[[a,b],b]}(-m)^* \right)$$

$$= \sum_{a,b} \left( \psi^{[a,b]}(-m) \wedge \psi^{[a,b]}(-m)^* - \psi^{[a,b]}(-m) \wedge \psi^{[a,b]}(-m)^* \right) = 0.$$

As the first factor represents the total co-adjoint action of $\mathfrak{g}$ on $\Lambda \otimes S$, $K \equiv 0$ on the subcomplex $C^\bullet$ of $\mathfrak{g}$-invariants, and our theorem is a special case of the following lemma. ✠

**(2.7.7) Lemma.** $\overline{\Box} = \Box + D + K$ *on* $\Lambda \otimes S$.

*Proof of the Lemma.* All the terms are second-order differential operators on $\Lambda \otimes S$. It suffices, then, to verify the identity on quadratic germs. The brutal calculations are performed in the Appendix. ✠

*Second proof of (2.7.5).* This is far cleaner, but relies on the Nakano identity proved in [T1]. Let *V* be a negatively graded $\mathfrak{g}[z]$-module, such that $z^m\mathfrak{g}$ maps $V(n)$ to $V(n+m)$. Assume that *V* carries a hermitian inner product, compatible with the hermitian involution on the zero-modes $\mathfrak{g} \subseteq \mathfrak{g}[z]$, and for which the graded pieces are mutually orthogonal. For us, *V* will be $S\mathfrak{g}[z]^*_{res}$. Write $R_a(m)$ for the action of $z^m\xi_a$ on *V* and define, for $m \geq 0$, $R_a(-m) = R_a(m)^*$. Define $\Box$ and $\overline{\Box}$ as before; our conditions on *V* ensure that they are well-defined. Define an endomorphism of $V \otimes \Lambda(\mathfrak{g}[z]/\mathfrak{g})^*_{res}$ by the formula

$$(2.7.8) \qquad T_V^\Lambda := \sum_{\substack{a,b \in A \\ m,n > 0}} \left\{[R_a(m), R_b(-n)] - R_{[a,b]}(m-n)\right\} \otimes \psi^a(-m) \wedge \psi^b(-n)^*.$$

Theorem (2.7.5) now splits up into the two propositions that follow; the first is known as *Nakano's Identity*, the second describes $T_V^\Lambda$ when $V = S\mathfrak{g}[z]^*_{res}$. ✠

**(2.7.9) Proposition.** *(cf. [T1], Prop. 2.4.7) On* $C^\bullet$, $\overline{\Box} = \Box + T_V^\Lambda + \deg$ *(*deg *is the exterior degree).* ✠

*(2.7.10) Remarks.* (i) Our $R_a(m)$ is the $\theta_a(m)$ of [T1], §2.4, whereas the operators $R_a(m)$ there are zero here (as is the level *h*). The constant $2c$ from [T1] is replaced by 1, because of our use of the Killing form, instead of the basic inner product. A sign discrepancy in the definition of $T_V^\Lambda$ arise because our $\xi_a$ here are self-adjoint, and not skew-adjoint as in [T1].
(ii) [T1] assumed finite-dimensionality of *V*, but the proof works equally well with our grading assumption.

**(2.7.11) Proposition.** *When* $V = S\mathfrak{g}[z]^*_{res}$, $T_V^\Lambda + \deg = D$ *on* $\Lambda \otimes S$.

*Proof.* Both sides are second-order differential operators on $\Lambda \otimes S$, so it suffices to check the equality on quadratic germs. Both sides annihilate $1 \otimes S$, so the following three checks suffice. In the first computation, we use the identity $\sum_a \psi^{[a,[a,b]]}(-n) = \psi^b(-n)$, coming from $\sum_a \operatorname{ad}(\xi_a)^2 = \mathbf{1}$ on $\mathfrak{g}$.

$$D\psi^b(-n) = \sum_{\substack{a \in A \\ 0 < m \leq n}} d_a(-m)\sigma^{[a,b]}(m-n)/n$$

$$= \sum_{\substack{a \in A \\ 0 < m \leq n}} \psi^{[a,[a,b]]}(m-n)/n = \psi^b(-n) = \left(T_V^\Lambda + \deg\right)\psi^b(-n);$$

$$D\left(\psi^b(-n) \wedge \psi^c(-p)\right) = D\psi^b(-n) \wedge \psi^c(-p) + \psi^b(-n) \wedge D\psi^c(-p)$$

$$= 2 \cdot \psi^b(-n) \wedge \psi^c(-p) = \left(T_S^\Lambda + \deg\right)\psi^b(-n) \wedge \psi^c(-p);$$

$$D\left(\sigma^c(-p) \cdot \psi^d(-q)\right) = \sigma^c(-p) \cdot D\psi^d(-q) + \frac{1}{q}\sum_{\substack{a \in A \\ 0 < m \leq q}} \sigma^{[a,d]}(m-q) \cdot \psi^{[a,c]}(-m-p),$$

which equals $\sigma^c(-p) \cdot \psi^d(-q)$ plus

$$T_S^\Lambda\left(\sigma^c(-p) \cdot \psi^d(-q)\right) = \frac{1}{q}\sum_{a \in A; 0 < m}\left\{[R_a(m), R_d(-q)] - R_{[a,d]}(m-q)\right\}\sigma^c(-p) \cdot \psi^a(-m)$$

$$= \frac{1}{q} \sum_{\substack{a \in A \\ 0 < m \leq p+q}} \sigma^{[a,[d,c]]}(m-q-p) \cdot \psi^a(-m) - \frac{1}{q} \sum_{\substack{a \in A \\ 0 < m \leq p}} \sigma^{[d,[a,c]]}(m-q-p) \cdot \psi^a(-m)$$

$$- \frac{1}{q} \sum_{\substack{a \in A \\ 0 < m \leq p+q}} \sigma^{[[a,d],c]}(m-q-p) \cdot \psi^a(-m)$$

$$= \frac{1}{q} \sum_{\substack{b \in A \\ 0 < m \leq p+q}} \sigma^{[d,[a,c]]}(m-q-p) \cdot \psi^a(-m) + \frac{1}{q} \sum_{\substack{a \in A \\ 0 < m \leq p}} \sigma^{[d,[a,c]]}(m-q-p) \cdot \psi^a(-m)$$

$$= \frac{1}{q} \sum_{\substack{a \in A \\ p < m \leq q}} \sigma^{[d,a]}(m-q-p) \cdot \psi^{[c,a]}(-m), \text{ as desired. } \maltese$$

*2.8. The harmonic forms and proof of (2.4.3)*

We now describe the harmonic forms, and therewith conclude the proof of Theorem (2.4.3). It will help to identify $\Lambda(\mathfrak{g}[z]/\mathfrak{g})^*_{res}$ with $\Lambda\mathfrak{g}[z]^*_{res}$, by using the isomorphism $d/dz : \mathfrak{g}[z]/\mathfrak{g} \cong \mathfrak{g}[z]$. This amounts to re-labeling the exterior generators, with $\psi^a(-m)$ now denoting what we used to call $(m+1) \cdot \psi^a(-m-1)$ ($m \geq 0$). Relations (2.6.1*) and (2.7.2*) now become:

(2.8.1)    $\mathrm{ad}_a(-m)^* \psi^b(-n) = \psi^{[a,b]}(m-n)$, or zero, if $m > n$,

(2.8.2)    $d_a(-m)^* \psi^b(-n) = \sigma^{[a,b]}(m-n-1)$, or zero, if $m > n+1$.

According to (2.7.5), the harmonic forms in the relative Koszul complex (2.6.4) are the forms in $\Lambda\mathfrak{g}[z]^*_{res} \otimes S\mathfrak{g}[z]^*_{res}$ killed by $d_a(-m)^*$, $m > 0$, and by $R_a(m) + \mathrm{ad}_a(-m)^*$, $m \geq 0$ ($a \in A$ in both cases).

Consider the *super-schemes* $\mathfrak{g}[[z,s]]$ and $G[[z,s]]$ of morphisms from the *formal super-disk* $\mathrm{Spf}\,\mathbb{C}[[z,s]]$ to $\mathfrak{g}$ and $G$, respectively. Concretely, $\mathfrak{g}[[z]]$ is the space of $\mathfrak{g}$-valued functions on the formal disk, and $\mathfrak{g}[[z,s]]$ is the graded vector space $\mathfrak{g}[[z]] \oplus s \cdot \mathfrak{g}[[z]]$. A polynomial on the latter involves only finitely many of the components $z^m\mathfrak{g}$, $sz^m\mathfrak{g}$. Similarly, $G[[z,s]] \cong G[[z]] \ltimes s\mathfrak{g}[[z]]$ as super-group scheme. We can identify our graded schemes with the tangent bundles to their even parts, by declaring the tangent spaces to be odd: thus, the *polynomials* on the super-schemes are identified with the *differential forms* on their even parts. (So the language of super-schemes is not strictly necessary in our case).

*(2.8.3) Remarks.* (i) For any scheme $X$, the regular functions on the super-scheme $X[s]$ of morphisms from the odd disk $\mathrm{Spec}\,\mathbb{C}[s]$ to $X$ are the algebraic differential forms on $X$.

(ii) Readers who prefer to think holomorphically may use the spaces $\mathfrak{g}\{z\}$ of germs of $\mathfrak{g}$-valued analytic maps. On holomorphic functions on $\mathfrak{g}\{z\}$, the $\mathbb{C}^\times$-scaling of $z$ normalizes the adjoint action of $G\{z\}$, and can be used to extract the invariant *polynomial* functions and forms.

**(2.8.4) Lemma.** *Identify $\Lambda\mathfrak{g}[z]^*_{res} \otimes S\mathfrak{g}[z]^*_{res}$ with the super-polynomials on $\mathfrak{g}[[z,s]]$. Then, the operators $d_a(-m)^*$, $m > 0$, and $R_a(m) + \mathrm{ad}_a(-m)^*$, $m \geq 0$, generate the co-adjoint action of $\mathfrak{g}[z,s]$.*

*Proof.* Clear from (2.8.1) and (2.8.2), as $d_a(-m)^*$ is the co-adjoint action of $s \cdot z^{m-1}\xi_a$. $\maltese$

**(2.8.5) Proposition.** *The harmonic forms in $C^\bullet$ are identified with the ring of polynomial functions on $\mathfrak{g}[[z,s]]$ which are invariant under the adjoint action of $G[[z,s]]$.*

*Proof.* The transition from $\mathfrak{g}[z,s]$-invariants to group invariants is obvious, because the co-adjoint action is locally-finite and factors, locally, through finite-dimensional truncations $\mathfrak{g}[z,s]/z^N$. $\maltese$

Call $P$ the *GIT quotient* $\mathfrak{g}//G := \mathrm{Spec}(S\mathfrak{g}^*)^G$; it can be identified with the space of primitives in the co-algebra $S\mathfrak{g}/[\mathfrak{g}, S\mathfrak{g}]$. The quotient map $q: \mathfrak{g} \to P$ induces a morphism $Q: \mathfrak{g}[[z,s]] \to P[[z,s]]$, which is invariant under the adjoint action of $G[[z,s]]$.

**(2.8.6) Proposition.** *The Ad-invariant polynomial functions on $\mathfrak{g}[[z,s]]$ are the pull-backs via Q of the polynomials on $P[[z,s]]$.*

*Proof.* Elements $\Lambda \mathfrak{g}[z]_{res}^* \otimes S\mathfrak{g}[z]_{res}^*$ are algebraic sections of the vector bundle $\Lambda \mathfrak{g}[z]_{res}^*$ over $\mathfrak{g}[[z]]$. As such, they are uniquely determined by their restriction to Zariski open subsets. The analogue holds for $P$. Now, the open subset $\mathfrak{g}^{rs} \subset \mathfrak{g}$ of *regular semi-simple* elements is an algebraic fibre bundle, via $q$, over the open subset $P^r \subset P$ of regular conjugacy classes. Let $\mathfrak{g}^{rs}[[z,s]]$ be the pull-back of $\mathfrak{g}^{rs}$ under the evaluation morphism $s = z = 0$. Because of the local product structure, it is clear that Ad-invariant polynomials over $\mathfrak{g}^{rs}[[z,s]]$ are precisely the $Q$-lifts of functions on $P^r[[z,s]]$. This shows, in particular, that the pull-back of polynomials from $P[[z,s]]$ to $\mathfrak{g}[[z,s]]$ is injective.

Let now $f$ be an invariant polynomial on $\mathfrak{g}[[z,s]]$. Its restriction to $\mathfrak{g}^{rs}[[z,s]]$ has the form $g \circ Q$, for some regular function $g$ on $P^r[[z,s]]$. Let $\mathfrak{g}^r \subset \mathfrak{g}$ be the open subset of *regular* elements. A theorem of Kostant's ensures that $q: \mathfrak{g}^r \to P$ is a submersion. In particular, it has local sections everywhere, so the morphism $Q: \mathfrak{g}^r[[z,s]] \to P[[z,s]]$ has local sections also. We can use local sections to extend our $g$ from $P^r[[z,s]]$ to $P[[z,s]]$, because $f$ was everywhere defined upstairs; the extension of $g$ is unique, and its lift via $Q$ must agree with $f$ everywhere, as it does so on an open set. So we have written $g$ as a pull-back. ✠

Here is an *algebraic* description of $Q$ on the *even* component $\mathfrak{g}[[z]]$. A polynomial $\Phi$ of degree $d$ on $P$ defines a pointwise evaluation map $P[[z]]^{\otimes d} \to \mathbb{C}[[z]]$. The $m$th coefficient $\Phi(-m)$ of the resulting power series is a polynomial on $P[[z]]$. Similarly, the polynomial $q^*\Phi := \Phi \circ q$ on $\mathfrak{g}$ gives a polynomial map from $\mathfrak{g}[[z]]$ to $\mathbb{C}[[z]]$, whose $m$th coefficient $(q^*\Phi)(-m)$ is a polynomial on $\mathfrak{g}[[z]]$. We then have $\Phi(m) \circ Q = (q^*\Phi)(m)$.

**(2.8.7) Proposition.** *Let $\Phi_1, \ldots, \Phi_\ell$ be a basis of linear functions on $P$ and let $\Phi_k(m)$ ($m \leq 0$) be the associated Fourier mode basis of linear functions on $P[[z]]$. After the $\psi$-relabeling, the cocycles $S_k(m)$ and $E_k(m)$ associated to $\Phi_k$ in (2.4.3) are the pull-backs under Q of $\Phi_k(m)$ and $d\Phi_k(m)$.*

*Proof.* We just checked this for $S_k(m)$; for $E_k(m)$, observe that the full $Q$ is the differential of its restriction to the even part, and that $E_k(m-1) = dS_k(m)$, after relabeling. ✠

Propositions (2.8.5), (2.8.6) and (2.8.7) together imply theorem (2.4.3).

*2.9. Cohomology of the super-Iwahori Lie algebra*
We now deduce Theorem (2.4.4) from (2.4.3). Let $\exp(\mathfrak{B})$ be the Iwahori subgroup of $G[[z]]$ associated to the $z$-adic completion $\mathfrak{B}_z$ of $\mathfrak{B}$. We shall write $H_{\mathcal{G}}^*(V)$ for *algebraic group cohomology* of $\mathcal{G}$ with coefficients in a representation $V$. Two applications of the van Est spectral sequence give

$$(2.9.1) \qquad H^*\left(\mathfrak{B}, \mathfrak{h}; S\mathfrak{B}_{res}^*\right) = H_{\exp(\mathfrak{B})}^*\left(S\mathfrak{B}_{res}^*\right), \qquad H_{res}^*\left(\mathfrak{g}[z], \mathfrak{g}; S\mathfrak{g}[z]_{res}^*\right) = H_{G[[z]]}^*\left(S\mathfrak{g}[z]_{res}^*\right).$$

The theorem follows now from *Shapiro's spectral sequence*

$$(2.9.2) \qquad E_2^{p,q} = H_{G[[z]]}^p\left(\mathbb{R}^q \mathrm{Ind}_{\exp(\mathfrak{B})}^{G[[z]]} S\mathfrak{B}_{res}^*\right) \Rightarrow H_{\exp(\mathfrak{B})}^{p+q}\left(S\mathfrak{B}_{res}^*\right),$$

in light of the following lemma, and the fact that $S\mathfrak{h}^*$ is a free $(S\mathfrak{g}^*)^{\mathfrak{g}}$-module.

**(2.9.3) Lemma.** $\mathbb{R}^q \mathrm{Ind}_{\exp(\mathfrak{B})}^{G[[z]]} S\mathfrak{B}_{res}^* = 0$ *for $q > 0$, while* $\mathrm{Ind}_{\exp(\mathfrak{B})}^{G[[z]]} S\mathfrak{B}_{res}^* = S\mathfrak{g}[z]_{res}^* \otimes_{(S\mathfrak{g}^*)^{\mathfrak{g}}} S\mathfrak{h}^*$, *with $G[[z]]$ acting on the right-hand side on the first factor only.*

*Proof.* $\mathbb{R}^q\mathrm{Ind}_{\exp(\mathfrak{B})}^{G[[z]]}\mathrm{S}\mathfrak{B}_{res}^*$ is the $q$th sheaf cohomology of the algebraic vector bundle $\mathrm{S}\mathfrak{B}_{res}^*$ over the quotient variety $G[[z]]/\exp(\mathfrak{B}) \cong G/B$, and also the $q$th cohomology of the structure sheaf $\mathcal{O}$ over the variety $G[[z]] \times_{Exp(\mathfrak{B})} \mathfrak{B}_z$, defined from the adjoint action of $\exp(\mathfrak{B})$ on $\mathfrak{B}_z$. This variety can be rewritten more transparently as $(G \times_B \mathfrak{b}) \times z\mathfrak{g}[[z]]$, whereby the left translation action of $G[[z]]$ becomes the $G$-action on the first factor (setting $z = 0$) coupled with conjugation on the second factor.

The first factor $G \times_B \mathfrak{b}$ maps properly and generically finitely to $\mathfrak{g}$, via $(g, \beta) \mapsto g\beta g^{-1}$. Using this, a theorem of Grauert and Riemenschneider guarantees the vanishing of the higher cohomology of $\mathcal{O}$. We identify the space of functions $G \times_B \mathfrak{b}$ with $\mathrm{S}\mathfrak{g}^* \otimes_{(\mathrm{S}\mathfrak{g}^*)^\mathfrak{g}} \mathrm{S}\mathfrak{h}^*$ thanks to the following Stein factorization of our map: $G \times_B \mathfrak{b} \to \mathfrak{h} \times_{\mathfrak{g}//G} \mathfrak{g} \to \mathfrak{g}$, where the first component of the first arrow is the projection $\mathfrak{b} \to \mathfrak{h}$, and second arrow projects onto the second factor. The left $G$-action on the variety becomes the adjoint action on $\mathrm{S}\mathfrak{g}^*$. It is now clear that the functions on $G[[z]] \times_{\exp(\mathfrak{B})} \mathfrak{B}_z$ are as stated in the Lemma. ✠

*(2.9.4) Remark.* We are really studying the conjugation morphism $G[[z]] \times_{\exp(\mathfrak{B})} \mathfrak{B}_z \to \mathfrak{g}[[z]]$ and its Stein factorization $G[[z]] \times_{\exp(\mathfrak{B})} \mathfrak{B}_z \to \mathfrak{h} \times_{\mathfrak{g}//G} \mathfrak{g}[[z]] \to \mathfrak{g}[[z]]$; Grauert-Riemenschneider implies higher cohomology vanishing. The construction is reduced to finite dimensions, by factoring out $z\mathfrak{g}[[z]]$ from the maps, and this justifies our argument.

**Appendix: Calculations for Lemma 2.7.7**

It is clear that both sides in (2.7.7) annihilate the constant line in $\Lambda \otimes S$, and it is also easy to see that they agree on the symmetric part $1 \otimes S$, where $D$, $K$, ad and ad$^*$ all vanish. So we must only check equality on the linear $\psi$ terms, and on the quadratic $\psi \wedge \psi$ and $\sigma\psi$ terms.

*A.1. The linear $\psi$ terms*

Fix $b \in A$, $n > 0$. We compute:

$$\bar{\partial}\psi^b(-n) = \frac{1}{2}\sum_{\substack{a \in A \\ 0 < m < n}} \psi^a(-m) \wedge \psi^{[a,b]}(m-n),$$

$$\bar{\partial}^*\bar{\partial}\psi^b(-n) = \frac{1}{4}\sum_{\substack{a \in A \\ 0 < m < n}} \frac{n}{m(n-m)}\psi^{[a,[a,b]]}(-n) - \frac{1}{4}\sum_{\substack{a \in A \\ 0 < m < n}} \frac{n}{m(n-m)}\psi^{[[a,b],a]}(-n) =$$

$$= \frac{1}{2}\sum_{\substack{a \in A \\ 0 < m < n}} \frac{n}{m(n-m)}\psi^{[a,[a,b]]}(-n) = \frac{1}{2}\sum_{0 < m < n} \frac{n}{m(n-m)}\psi^b(-n) =$$

$$= \frac{1}{2}\sum_{0 < m < n}\left(\frac{1}{m} + \frac{1}{n-m}\right)\cdot \psi^b(-n) = \sum_{0 < m < n} 1/m \cdot \psi^b(-n)$$

Further, $\bar{\partial}^*\psi^b(-n) = 0$, so $\square\psi^b(-n)$ is as just computed. Next,

$$\square\psi^b(-n) = \sum_{\substack{a \in A \\ 0 < m < n}} \frac{1}{m}\mathrm{ad}_a(-m)\mathrm{ad}_a(-m)^*\psi^b(-n) = \sum_{\substack{a \in A \\ 0 < m < n}} \frac{1}{m}\cdot\frac{n-m}{n}\cdot\psi^{[a,[a,b]]}(-n) =$$

$$= \sum_{0 < m < n}\left(\frac{1}{m} - \frac{1}{n}\right)\cdot\psi^b(-n) = \sum_{0 < m < n} 1/m \cdot \psi^b(-n) - \psi^b(-n) + \psi^b(-n)/n$$

$$D\psi^b(-n) = \sum_{\substack{a \in A \\ 0 < m < n}} \psi^{[a,[a,b]]}(-n)/n = \psi^b(-n),$$

$$K\psi^b(-n) = \frac{1}{n}\sum_{a\in A} \mathrm{ad}_{[a,b]}(0)\psi^a(-n) = -\frac{1}{n}\sum_{a\in A} \psi^{[a,[a,b]]}(-n) = -\psi^b(-n)/n,$$

and the last three terms sum up to $\overline{\Box}\psi^b(-n)$, as claimed. ✠

*A.2. The quadratic $\psi \wedge \psi$ terms*

Fix $b, c \in A$ and $n, p > 0$. For each second-order differential operator $P$ involved, we focus on the *cross-term* $P(\psi^b(-n) \wedge \psi^c(-p)) - P\psi^b(-n) \wedge \psi^c(-p) - \psi^b(-n) \wedge P\psi^c(-p)$; equality of cross-terms and the identity (2.7.7) on the linear factors imply the identity for quadratic terms.

**(A.3) Claim.** *The cross-term in $\overline{\Box}(\psi^b(-n) \wedge \psi^c(-p))$ is the following sum of three terms:*

$$\sum_{\substack{a\in A \\ 0<m<n}} \left(\frac{1}{m} - \frac{1}{n}\right) \cdot \psi^{[a,b]}(m-n) \wedge \psi^{[a,c]}(-m-p) + \sum_{\substack{a\in A \\ 0<m<p}} \left(\frac{1}{m} - \frac{1}{p}\right) \cdot \psi^{[a,b]}(-m-n) \wedge \psi^{[a,c]}(m-p)$$

$$-\sum_{a\in A} \left(\frac{1}{n} + \frac{1}{p}\right) \cdot \psi^{[a,b]}(-n) \wedge \psi^{[a,c]}(-p).$$

*Proof.* We rewrite the sum, by adding and subtracting terms:

(A.4)
$$\sum_{\substack{a\in A \\ 0<m<n}} \left(\frac{1}{m} + \frac{1}{p}\right) \cdot \psi^{[a,b]}(m-n) \wedge \psi^{[a,c]}(-m-p) + \sum_{\substack{a\in A \\ 0<m<p}} \left(\frac{1}{m} + \frac{1}{n}\right) \cdot \psi^{[a,b]}(-m-n) \wedge \psi^{[a,c]}(m-p)$$

$$- \sum_{\substack{a\in A \\ 0<m<n+p}} \left(\frac{1}{n} + \frac{1}{p}\right) \cdot \psi^{[a,b]}(-m) \wedge \psi^{[a,c]}(m-n-p).$$

We now track, in turn, the source of each of the three terms in (A.4). We have

(A.5)
$$\overline{\partial}(\psi^b(-n) \wedge \psi^c(-p)) = \frac{1}{2} \sum_{\substack{a\in A \\ 0<m<n}} \psi^a(-m) \wedge \psi^{[a,b]}(m-n) \wedge \psi^c(-p)$$

$$+ \frac{1}{2} \sum_{\substack{a\in A \\ 0<m<p}} \psi^a(-m) \wedge \psi^b(-n) \wedge \psi^{[a,c]}(m-p),$$

and applying $\overline{\partial}^*$ to the first term gives

(A.6)
$$\overline{\partial}^*\overline{\partial}\psi^b(-n) \wedge \psi^c(-p)$$
$$+ \frac{1}{4} \sum_{\substack{a\in A \\ 0<m<n}} \frac{1}{m}\psi^{[a,b]}(m-n) \wedge \mathrm{ad}_a(m)^*\psi^c(-p) - \frac{1}{4} \sum_{\substack{a\in A \\ 0<m<n}} \frac{1}{n-m}\psi^a(-m) \wedge \mathrm{ad}_{[a,b]}(n-m)^*\psi^c(-p)$$
$$+ \frac{1}{4} \sum_{\substack{a\in A \\ 0<m<n}} \frac{1}{p}\mathrm{ad}_c(p)^*(\psi^a(-m) \wedge \psi^{[a,b]}(m-n)).$$

The first term is $\overline{\Box}\psi^b(-n) \wedge \psi^c(-p)$. The second-line terms agree, after substituting $m \leftrightarrow (n-m)$, and sum to

(A.7)
$$\frac{1}{2} \sum_{\substack{a\in A \\ 0<m<n}} \frac{p+m}{mp} \psi^{[a,b]}(m-n) \wedge \psi^{[a,c]}(-m-p).$$

Amusingly, the third line takes the same value (A.7). So the sum in (A.6) equals

(A.8) $$\Box\psi^b(-n)\wedge\psi^c(-p) + \sum_{\substack{a\in A \\ 0<m<n}} \frac{p+m}{mp}\psi^{[a,b]}(-m-p)\wedge\psi^{[a,c]}(m-n),$$

and so the cross-term in (A.8) accounts for the first term in (A.4).

Substituting $b\leftrightarrow c$, $n\leftrightarrow p$ shows that the $\bar\partial^*$-image of the second term in (A.5) is

(A.9) $$\psi^b(-n)\wedge\Box\psi^c(-p) + \sum_{\substack{a\in A \\ 0<m<n}} \frac{n+m}{mn}\psi^{[a,b]}(-m-n)\wedge\psi^{[a,c]}(m-p),$$

whose cross-term is the second term in (A.4).

Finally, $\bar\partial^*\bigl(\psi^b(-n)\wedge\psi^c(-p)\bigr) = \frac{n+p}{np}\cdot\psi^{[b,c]}(-p-n)$, whence by applying $\bar\partial$ we get

$$\bar\partial\bar\partial^*\bigl(\psi^b(-n)\wedge\psi^c(-p)\bigr) = \frac{1}{2}\frac{n+p}{np}\sum_{\substack{a\in A \\ 0<m<n+p}} \psi^a(-m)\wedge\psi^{[a,[b,c]]}(m-p-n) =$$

$$= \frac{n+p}{2np}\sum_{\substack{a\in A \\ 0<m<n+p}} \bigl(\psi^a(-m)\wedge\psi^{[[a,b],c]}(m-p-n) + \psi^a(-m)\wedge\psi^{[b,[a,c]]}(m-p-n)\bigr) =$$

$$= -\frac{n+p}{np}\sum_{\substack{a\in A \\ 0<m<n+p}} \psi^{[a,b]}(-m)\wedge\psi^{[a,c]}(m-p-n)$$

which is the third term in (A.4). The claim is proved. ✠

Now $D\bigl(\psi^b(-n)\wedge\psi^c(-p)\bigr) = D\psi^b(-n)\wedge\psi^c(-p) + \psi^b(-n)\wedge D\psi^c(-p)$, with no cross-term; while

(A.10) $$K\bigl(\psi^b(-n)\wedge\psi^c(-p)\bigr) = K\psi^b(-n)\wedge\psi^c(-p) + \psi^b(-n)\wedge K\psi^c(-p)$$
$$-\frac{1}{n}\psi^{[a,b]}(-n)\wedge\psi^{[a,c]}(-p) - \frac{1}{p}\psi^{[a,b]}(-n)\wedge\psi^{[a,c]}(-p),$$

(A.11) $$\Box\bigl(\psi^b(-n)\wedge\psi^c(-p)\bigr) = \Box\psi^b(-n)\wedge\psi^c(-p) + \psi^b(-n)\wedge\Box\psi^c(-p)$$
$$+ \sum_{\substack{a\in A \\ 0<m<n}} \frac{1}{m}\frac{n-m}{n}\psi^{[a,b]}(m-n)\wedge\psi^{[a,c]}(-m-p)$$
$$+ \sum_{\substack{a\in A \\ 0<m<p}} \frac{1}{m}\frac{p-m}{p}\psi^{[a,b]}(-m-n)\wedge\psi^{[a,c]}(m-p),$$

and the cross-terms in (A.10) and (A.11) add up to the expression in (A.3). ✠

*A.12. The $\sigma\psi$ terms*
As before, fix $b,c\in A$ and $n,p>0$. Then,

$$\bar\partial\bigl(\sigma^b(-n)\psi^c(-p)\bigr) = \sum_{\substack{a\in A \\ 0<m\le n}} \sigma^{[a,b]}(m-n)\cdot\psi^a(-m)\wedge\psi^c(-p) + \frac{1}{2}\sum_{\substack{a\in A \\ 0<m<p}} \sigma^b(-n)\psi^a(-m)\wedge\psi^{[a,c]}(m-p),$$

and applying $\bar\partial^*$ yields the following sum:

$$\overline{\Box}\sigma^b(-n)\cdot\psi^c(-p) - \frac{1}{p}\sum_{\substack{a\in A\\0<m\leq n}}\sigma^{[a,c]}(m-p-n)\psi^{[a,b]}(-m)$$

$$+\frac{1}{2}\sum_{\substack{a\in A\\0<m<p}}\frac{1}{m}\sigma^{[a,b]}(-m-n)\psi^{[a,c]}(m-p) - \frac{1}{2}\sum_{\substack{a\in A\\0<m<p}}\frac{1}{m-p}\sigma^{[[a,c],b]}(m-p-n)\psi^a(-m)$$

$$+\frac{1}{2}\sum_{\substack{a\in A\\0<m\leq n}}\frac{m+p}{mp}\sigma^{[a,b]}(m-n)\psi^{[a,c]}(-m-p) + \frac{1}{2}\sum_{\substack{a\in A\\0<m\leq n}}\frac{m+p}{mp}\sigma^{[a,b]}(m-n)\psi^{[a,c]}(-m-p)$$

$$+\sigma^b(-n)\overline{\Box}\psi^c(-n).$$

The first two lines come from the $R^*$-terms in $\overline{\partial}^*$, the last two lines from the $\mathrm{ad}^*$-terms. The two terms in each of the middle rows are equal, so the cross-term can be rewritten as follows:

(A.13)
$$-\frac{1}{p}\sum_{\substack{a\in A\\0<m\leq n}}\sigma^{[a,c]}(m-p-n)\cdot\psi^{[a,b]}(-m) + \sum_{\substack{a\in A\\0<m<p}}\frac{1}{m}\sigma^{[a,b]}(-m-n)\cdot\psi^{[a,c]}(m-p)$$

$$+\sum_{\substack{a\in A\\0<m\leq n}}\left(\frac{1}{m}+\frac{1}{p}\right)\sigma^{[a,b]}(m-n)\cdot\psi^{[a,c]}(-m-p)$$

Now, $\overline{\partial}^*\bigl(\sigma^b(-n)\cdot\psi^c(-p)\bigr) = \sigma^{[c,b]}(-n-p)/p$, whence

(A.14)
$$\overline{\partial}\,\overline{\partial}^*\bigl(\sigma^b(-n)\cdot\psi^c(-p)\bigr) = \frac{1}{p}\sum_{\substack{a\in A\\0<m\leq n+p}}\sigma^{[a,[c,b]]}(m-n-p)\cdot\psi^a(-m)$$

$$= \frac{1}{p}\sum_{\substack{a\in A\\0<m<n+p}}\sigma^{[a,c]}(m-n-p)\cdot\psi^{[a,b]}(-m)$$

$$-\frac{1}{p}\sum_{\substack{a\in A\\0<m<n+p}}\sigma^{[a,b]}(m-n-p)\cdot\psi^{[a,c]}(-m).$$

Summing (A.13) and (A.14) gives

(A.15)
$$\frac{1}{p}\sum_{\substack{a\in A\\0<m\leq p}}\sigma^{[a,c]}(m-p)\cdot\psi^{[a,b]}(-m-n) + \sum_{\substack{a\in A\\0<m<p}}\left(\frac{1}{m}-\frac{1}{p}\right)\sigma^{[a,b]}(-m-n)\cdot\psi^{[a,c]}(m-p)$$

$$+\sum_{\substack{a\in A\\0<m<n}}\frac{1}{m}\sigma^{[a,b]}(m-n)\cdot\psi^{[a,c]}(-m-p) - \frac{1}{p}\sum_{a\in A}\sigma^{[a,b]}(-n)\cdot\psi^{[a,c]}(-p);$$

here, the first term is the sum of the first terms in (A.13) and (A.14), the second and third incorporate the second and third terms in (A.13) and the $0 < m < p$, resp. the $p < m < p+n$ portions of the second term in (A.14), and the final term is the $m = p$ contribution of the same.

Moving on to the right-hand side of (2.7.7), the cross-term in $\Box\bigl(\sigma^b(-n)\psi^c(-p)\bigr)$ is

(A.16)
$$\sum_{\substack{a\in A\\0<m<n}}\frac{1}{m}\sigma^{[a,b]}(m-n)\cdot\psi^{[a,c]}(-m-p) + \sum_{\substack{a\in A\\0<m<p}}\frac{p-m}{mp}\sigma^{[a,b]}(-m-n)\cdot\psi^{[a,c]}(m-p),$$

the two terms coming from the $\mathrm{ad}\cdot R^*$ and $R\cdot\mathrm{ad}^*$ cross-terms, respectively. Further,

$$\text{(A.17)} \quad D\left(\sigma^b(-n)\cdot\psi^c(-p)\right) = \sigma^b(-n)\cdot D\psi^c(-p) + \frac{1}{p}\sum_{\substack{a\in A \\ 0<m\leq p}} \sigma^{[a,c]}(m-p)\cdot\psi^{[a,b]}(-m-n),$$

$$\text{(A.18)} \quad K\left(\sigma^b(-n)\cdot\psi^c(-p)\right) = \sigma^b(-n)\cdot K\psi^c(-p) - \frac{1}{p}\sum_{a\in A} \sigma^{[a,b]}(-n)\cdot\psi^{[a,c]}(-p).$$

It is now clear that the cross-terms in (A.16)–(A.18) sum to (A.15). ✠

## References

[F1] Feigin, B.I.: On the cohomology of the Lie algebra of vector fields and of the current algebra. Sel. math. Sov. **7** (1988), 49–62

[F2] Feigin, B.I.: Differential operators on the moduli space of G-bundles over curves and Lie algebra cohomology. Special functions (Okayama, 1990), ICM-90 Satell. Conf. Proc., Springer, Tokyo, (1991), 90–103

[FGZ] Frenkel, I.B., Garland, H., and Zuckerman, G.J.: Semi-infinite cohomology and string theory. Proc. Nat. Acad. Sci. U.S.A. **83** (1986), 8442–8446

[H] Hanlon, P.: Some conjectures and results concerning the homology of nilpotent Lie algebras. Adv. Math. **84** (1990), 91–134

[K] Koszul, J-.L.: Sur la cohomologie relative des algèbres de Lie. C. R. Acad. Sci. Paris **228**, (1949), 457–459

[M] Macdonald, I.G.: A formal identity for affine root systems. Unpublished manuscript

[T1] Teleman, C.: Lie algebra Cohomology and the Fusion Rules, Commun. Math. Phys. **173** (1995), 265–311

[T2] Teleman, C.: Some Hodge theory from Lie algebras. To appear in: Proceedings of the UCI conference on Motives, Polylogarithms and Non-Abelian Hodge Theory, June 1998.

Susanna Fishel  
Mental Images  
Fasanenstraße 81  
10623 Berlin, Germany  

Ian Grojnowski  
DPMMS  
Centre for Mathematical Sciences  
Wilberforce Road  
Cambridge, CB3 0WB, UK  

Constantin Teleman  
DPMMS  
Centre for Mathematical Sciences  
Wilberforce Road  
Cambridge, CB3 0WB, UK